\documentclass[11pt]{article}

\usepackage[a4paper, margin=2.5cm]{geometry}

\usepackage[T1]{fontenc}
\usepackage{lmodern}
\usepackage{microtype}

\usepackage{amsmath, amssymb, mathtools}

\usepackage{graphicx}
\usepackage{xcolor}
\usepackage{booktabs}
\usepackage{subcaption}
\usepackage{tabularx}
\usepackage{tikz}
\usepackage{enumitem}
\usepackage{siunitx}
\usepackage[numbers]{natbib}

\usepackage{algorithm}
\usepackage[noend]{algpseudocode}

\usepackage[colorlinks=true,linkcolor=blue,citecolor=blue,urlcolor=blue]{hyperref}
\usepackage[nameinlink,capitalize,noabbrev]{cleveref}

\crefname{equation}{Eq.}{Eqs.}
\crefname{figure}{Fig.}{Figs.}
\crefname{table}{Table}{Tables}
\crefname{algorithm}{Algorithm}{Algorithms}
\crefname{section}{Section}{Sections}

\makeatletter
\@ifpackageloaded{algpseudocode}{%
    \providecommand{\STATE}{\State}%
    \providecommand{\IF}{\If}%
    \providecommand{\ENDIF}{\EndIf}%
    \providecommand{\ELSE}{\Else}%
}{}
\makeatother

\newcommand{\feelpp}{\texttt{Feel++}}

\newcommand{\Solid}{\mathcal{S}} 
\newcommand{\Fluid}{\mathcal{F}} 
\newcommand{\Density}{\rho} 
\newcommand{\Real}{\mathbb{R}} 
\newcommand{\Alemap}{\mathcal{A}} 
 
\newcommand{\ALE}{\mathcal{A}} 

\title{Two-way Coupling of Fluid--Structure Interaction
for Elastic Magneto-Swimmers:
A Finite Element ALE Approach}

\author{
  Christophe Prud'homme\thanks{IRMA, University of Strasbourg \& CEMOSIS, France.} \\
  \and Vincent Chabannes\footnotemark[1] \\
  \and Laetitia Giraldi\thanks{Inria Sophia Antipolis Méditerranée, Calisto Team, France.} \\
  \and Agathe Chouippe\thanks{Université de Strasbourg, CNRS ICUBE UMR 7357, France.} \\
  \and Céline Van Landeghem\footnotemark[1]
}

\date{} 

\begin{document}

\maketitle

\begin{abstract}
Artificial micro-swimmers actuated by external magnetic fields hold significant promise for targeted biomedical applications, including drug delivery and micro-robot-assisted therapy. However, their dynamics remain challenging to control due to the complex nonlinear coupling between magnetic actuation, elastic deformations, and fluid interactions in confined biological environments. Numerical modeling is therefore essential to better understand, predict, and optimize their behavior for practical applications.
In this work, we present a comprehensive finite element framework based on the Arbitrary Lagrangian–Eulerian formulation to simulate deformable elastic micro-swimmers in confined fluid domains. The method employs a full-order model that resolves the complete fluid dynamics while simultaneously tracking swimmer deformation and global displacement on conforming meshes. 
Numerical experiments are performed with the open-source finite element library $\feelpp$, demonstrating excellent agreement with experimental data from the literature. The validation benchmarks in both two and three dimensions confirm the accuracy, robustness, and computational efficiency of the proposed framework, representing a foundational step toward developing digital twins of magneto-swimmers for biomedical applications.
\end{abstract}


\section{Introduction}
\label{introduction}


In recent years, significant attention has been devoted to the development of 
artificial micro-swimmers for biomedical applications, such as targeted drug 
delivery or minimally invasive diagnostics \cite{bunea_recent_2020}. One promising 
approach involves actuating these micro-swimmers using external magnetic fields, 
often generated by MRI technology \cite{feng_mini_2014}, to enable wireless control 
inside the human body. However, controlling flagellated swimmers in such conditions 
remains a challenging task. These systems are typically non-homologous, and their 
dynamics are often poorly controllable \cite{giraldi_local_2016,alouges_purcells_2017} 
due to the complex coupling between magnetic actuation, elastic deformations, 
and hydrodynamic interactions at low Reynolds number. This complexity motivates 
the development of accurate numerical models capable of capturing the coupled 
fluid-elasto-magneto interactions governing their dynamics.

Modeling the motion of micro-swimmers in fluids involves the coupling of fluid 
dynamics with rigid and, eventually, elastic body mechanics, as well as complex 
fluid–structure interactions. Various numerical approaches have been developed 
to address this problem, each offering a trade-off between computational efficiency 
and modeling accuracy. 

At low Reynolds numbers, where inertial effects are negligible, simplified models, 
such as the \textit{Resistive Force Theory} 
\cite{gray_propulsion_1955,el_alaoui-faris_optimal_2020,jikeli_sperm_2015} 
and the \textit{Slender Body Theory} 
\cite{cox_motion_1970, batchelor_slender-body_1970, keller_slender-body_1976, johnson_improved_1980}, 
provide computationally efficient approximations. These methods are particularly useful 
for slender filaments, but they rely on strong geometric assumptions and do not fully 
capture near-field hydrodynamic interactions or swimmer–boundary effects.

To improve accuracy, \textit{Boundary Integral Methods} 
\cite{shum_microswimmer_2019,walker_boundary_2019} have been introduced, offering a 
more refined description of fluid-structure coupling without discretization of the 
entire fluid domain \cite{pozrikidis_boundary_1992}. These methods rely on Green's 
functions of the Stokes equations \cite{blake_note_1971}, which are available in 
closed form and allow the velocity of the fluid to be expressed as an integral 
over the swimmer's surface. However, this approach is restricted to the Stokes 
regime and cannot be generalized to inertial flows. In addition, the singular 
behavior of the Green functions near the boundary requires regularization 
techniques to ensure numerical stability and precision 
\cite{lefebvre-lepot_accurate_2015,huang_notes_1993,olson_coupling_2011}.

When more complex geometries, large deformations, or interactions with boundaries 
are involved, full discretization of both fluid and swimmer domains becomes necessary. 
In such cases, the \textit{Finite Element Method} provides a versatile framework, 
particularly when combined with techniques like the \textit{Immersed Boundary Method }
\cite{bergmann_modeling_2011,morab_overview_2020,kim_immersed_2019,bergmann_accurate_2014, bergmann_bioinspired_2016} 
or \textit{CutFEM} \cite{monasse_conservative_2012, hansbo_cut_2016, burman_cutfem_2015}. 
These approaches allow the simulation of moving rigid swimmers immersed in a fluid, but they 
introduce additional challenges, such as force interpolation errors and numerical 
instabilities near the fluid–structure interface. Overall, while existing models 
provide valuable insight into micro-swimmer dynamics, many approaches that include 
both fluid- and elastic-solid interactions are limited to configurations where the 
solid undergoes only passive deformations, without rigid-body motion. As a result, 
these models do not capture the full dynamics of a self-propelled swimmer, 
particularly in scenarios where elastic deformation allows the swimmer to push 
the fluid, enabling its displacement.

In this work, we introduce a finite element method based on the 
\textit{Arbitrary Lagrangian–Eulerian} framework to compute the displacement of 
deformable swimmers immersed in viscous flows within geometrically complex domains 
\cite{chabannes_high-order_2013, pena_high-order_2012,berti_numerical_2021,van_landeghem_mathematical_2024,van_landeghem_towards_2024}. 
This approach resolves the full fluid dynamics while tracking both swimmer 
deformation and global motion on a conforming mesh. Compared to immersed boundary 
methods, the explicit tracking of the interface during the computation allows for 
precise evaluation of all physical quantities at the fluid–structure boundary. 
Thus, it provides an accurate description of the fluid–structure interface and 
can be extended to complex fluids relevant to biological applications.

The paper is organized as follows. In Section \ref{modeling} we introduce the 
mathematical models, starting with the fluid equations, then describing the 
swimmer through rigid-body motion and elasticity, and finally presenting the 
coupled fluid–structure interaction problem. Section \ref{discretization} 
details the spatial and temporal discretization of all these equations. 
Section \ref{strategies} focuses on computational strategies, including the 
remeshing procedure, and the preconditioning technique. Section \ref{results} 
presents the numerical results, consisting of validation tests, all performed 
using the open-source finite element library Feel++ \cite{feel_team_feel_nodate}. 
Finally, section \ref{conclusion} provides conclusions and perspectives.



\section{Mathematical modeling}
\label{modeling}

This section introduces the mathematical model for the fluid–magneto-swimmer interaction problem. We first describe the hydrodynamics of incompressible Newtonian fluids and the dynamics of elastic bodies, governed by hyper-elasticity equations. 
We then present the coupling between these two components, taking into account the rigid motion of the swimmer.

\subsection{Fluid model}

Let $\Fluid^t \subset \Real^d$, with $d=2,3$ denoting the spatial dimension, be the region 
occupied by the fluid at time $t \in \, ]0,T]$, where $T > 0$ is the final time. 
In what follows, we use the notation with superscript $t$ to emphasize that the corresponding quantity is time-dependent.
The physical properties of the incompressible Newtonian fluid are its dynamic viscosity 
$\mu_{\Fluid} \in \Real^+$ and its density $\Density_{\Fluid} \in \Real^+$.
Their hydrodynamics are described by the non-linear Navier-Stokes equations. 
These equations are formulated in the Eulerian 
framework and describe the time evolution of the velocity field 
$u^t: \, ]0,T]\times \Fluid^t \rightarrow \Real^d$, 
and the pressure field $p^t: \, ]0,T]\times\Fluid^t \rightarrow \Real$ of the fluid.

The Navier-Stokes equations are given at time $t > 0$ by
\small
\begin{equation*}
	\left\{
	\begin{aligned}
		\rho_\Fluid \Big(\partial_t u^t + ( u^t \cdot \nabla )u^t \Big) - \nabla \cdot \sigma(u^t,p^t) &=  0_{\Real^d}
		\, && \text{in $ \Fluid^t$},\\
		\nabla \cdot u^t &= 0  \, && \text{in $ \Fluid^t$},\\
        u^t &= 0_{\Real^d} \quad && \text{on $\partial \Fluid^t_D$},\\
		\sigma(u^t,p^t) n_{\Fluid}^t &=  0_{\Real^d} && \text{on $ \partial \Fluid^t_N$}.
	\end{aligned}
	\right.
\end{equation*}
\normalsize
The boundary conditions are either homogeneous Dirichlet boundary conditions, 
imposed on $\partial \Fluid^t_D$, or homogeneous Neumann boundary conditions, 
imposed on $\partial \Fluid^t_N$, such that $\partial \Fluid^t = \partial \Fluid^t_D \cup \partial \Fluid^t_N$. 
Here, $n_{\Fluid}^t$ is the outward unit normal vector to $\partial \Fluid^t_N$. 
The stress tensor is defined as
\begin{align*}
	\sigma(u^t,p^t) &= -p^t \mathbb{I}_d +  \mu_{\Fluid} \bigl( \nabla u^t + (\nabla u^t)^T \bigr),
\end{align*}
where $\mathbb{I}_d$ is the identity matrix of size $d$.

\subsection{Elasto-magneto-swimmer model}

\begin{figure}[!ht]
\begin{center}
\begin{tikzpicture}[scale = 2]

  \draw[thick, black] (-1.5,-1.2) rectangle (4.8,2);
  \node[anchor=south west] at (-1.5,-1.2) {\small$\Fluid^t$};

  \begin{scope}[rotate=15]
    \draw[color=black, ->, thick] (-0.8,0.5) -- (-1.2,0.5);
  \end{scope}
  \node[color=black] at (-1.05, 0.) {\small{$B^t$}};

  \begin{scope}[rotate=15]
    \filldraw[fill=gray, fill opacity=0.2, draw=black, thick] (0,0) rectangle (0.5,1);
    \draw[black, thick, dotted] (0.5,0) -- (0.5,1);
    \node at (0.25,0.25) {\color{black}\small{$\Solid_{\text{head}}^t$}};
    \node at (1.5,0.3) {\color{black}\small{$\Solid_{\text{tail}}^t$}};

    \draw[color=black, thick] (0.5,1) .. controls (2.2,0.4) .. (4,-1);
    \draw[color=black, thick] (0.5,0) .. controls (2.2,-0.1) .. (4,-1);
    \node[color=black] at (2, -0.35) {\small{$\Solid^t$}};

    \draw[->, color=black, thick] (0,0.5) -- (-0.6,0.5);
    \node at (-0.2, 0.35) {\color{black}\small{$M^t$}};
  \end{scope}

  \draw[dotted, thick] (-0.6,0.5) -- (0.6,0.5);

  \draw[black,thick] (-0.5,0.5) arc[start angle=180,end angle=245,radius=0.1];
  \node at (-0.4,0.7) {\color{black}\small$\theta^t$};

\end{tikzpicture}
\end{center}

\caption{Configuration of the magneto-swimmer at time $t$ inside a fluid domain 
$\Fluid^t$. The magnetic moment $M^t$ tends to align with the external magnetic field $B^t$, 
and $\theta^t$ represents the swimmer’s head orientation.}
\label{fig:magneto-swimmer_2}
\end{figure}
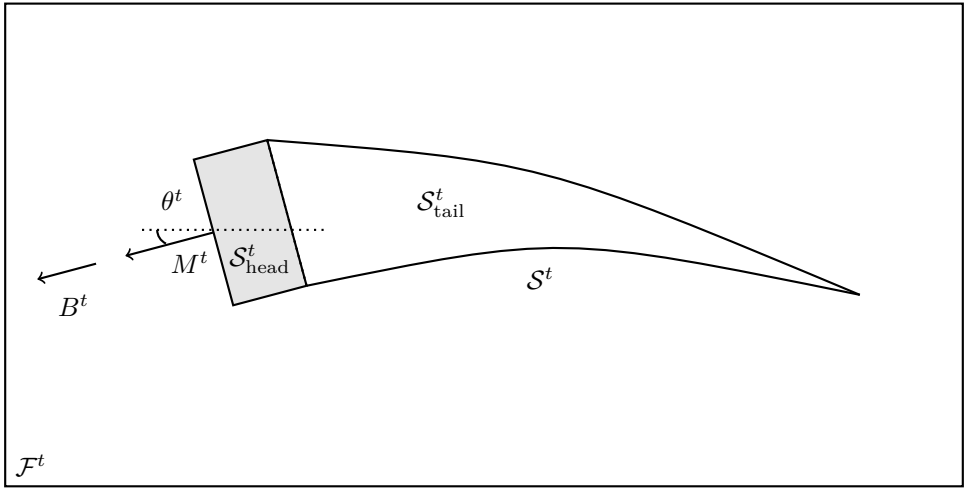

The swimmer \(\Solid\) consists of its magnetic head 
\(\Solid_{\text{head}}\) and its elastic tail \(\Solid_{\text{tail}}\), such that 
\(\Solid = \Solid_{\text{head}} \cup \Solid_{\text{tail}}\). 
We introduce the following notation: $\Solid^*$ denotes the reference domain of the 
swimmer, i.e., the initial domain at each remeshing step, while $\Solid^t$ denotes 
the domain at the current time. The same notation is used for the fluid. 
The notations are shown in \cref{fig:magneto-swimmer_2}.
The swimmer is described by the position of its center of mass $x^t_{cm}$ and the 
configuration of its structure, i.e., its deformation. These two quantities are 
determined by its time-dependent translational velocity $U^t$, 
 angular velocity $\omega^t$, and the elastic displacement $\eta^t$ of the flagella. 
 This elastic displacement is actuated by an external magnetic field 
$B^t : \, ]0,T] \to \Real^d$.\\

The linear velocity $U^t : \, ]0,T] \to \Real^d$ is described by the Newton's second law 
\begin{equation}
	\begin{aligned}
		m^t \, \mathrm{d}_{t} U^t &=  f^t_H,\\
	\end{aligned}
	\label{Eq:translation}
\end{equation}

\medskip 

\noindent and the angular velocity $\omega^t : \, ]0,T] \to \Real^{d^*}$, 
with $d^* = 1$ in two dimensions and $d^* = 3$ in three dimensions,
 is described by the Euler's equation
\begin{equation}
	\begin{aligned}
		\mathrm{d}_{t} [R(\theta^t) J^t R(\theta^t)^T \omega^t] &= T_H ^t + T^t_m.
	\end{aligned}
	\label{Eq:angular}
\end{equation}

\medskip 

\noindent In \eqref{Eq:translation}, $m^t \in \Real^+$ describes the mass of the swimmer, 
which is computed from its density $\rho_{\Solid} \in \Real^+$.
In \eqref{Eq:angular}, $J^t \in S_{++}^{d^*}$ represents its inertia tensor,
$R$ the rotation matrix, and $\theta^t : \, ]0,T] \to \Theta$ 
the rotation 
angle of the swimmer, which is derived from
$$
\mathrm{d}_{t} \theta^t = \omega^t.
$$
In two dimensions, the interval $\Theta = [-\pi,\pi]$, whereas in three dimensions 
it is given by $\Theta = [-\pi, \pi] \times [0, \pi] \times [0, \pi/2]$.
In \eqref{Eq:translation} and \eqref{Eq:angular}, $f^t_H :  \, ]0,T] \times \partial \Solid^t \to \Real^d$ and $T^t_H :  \, ]0,T] \times \partial \Solid^t \to \Real^{d^*}$ describe the hydrodynamical forces and torques
\begin{align}
\label{eq:hydro_forces}
f^t_{H} &= \int_{\partial \Solid^t} \sigma(u^t,p^t) n^t_{\Solid}, \quad
T^t_H = \int_{\partial \Solid^t} (x^t - x^t_{cm}) \times (\sigma(u^t,p^t) n^t_{\Solid}),
\end{align}
where $n^t_{\Solid}$ is the current outward unit normal vector to $\partial \Solid^t$.
The swimmer is also subjected to the magnetic torque $T_m^t :  \, ]0,T] \times \Solid^*_{\text{head}} \to \Real^{d^*}$ defined as
$$
T_m^t = \mathsf{m} R(\theta^t) M^t \times B^t,
$$
where $M^t :  \, ]0,T] \to \Real^d$ is its magnetic moment and $\mathsf{m} \in \Real^+$ its magnetization.

Moreover, the elastic displacement $\eta^t :  \, ]0,T] \times \Solid^* \to \Real^d$ of 
the flagella is governed by the hyper-elasticity equations
\small
\begin{equation}
\label{eq:hyper-elasticity}
	\left\{
	\begin{aligned}
		\rho_{\Solid} \, \partial_{tt} \eta^t - \nabla \cdot \big( F(\eta^t) \, \Sigma(\eta^t) \big) &= 0_{\Real^d}  
		 && \text{in $ \Solid^*$},\\
		\eta^t &= R(\theta^t) (x^* - x^*_{cm}) - (x^* - x^*_{cm}) 
		 && \text{on $ \Solid^*_{\text{head}}$},\\
		 F(\eta^t) \, \Sigma(\eta^t)  n^*_{\Solid} &= 0_{\Real^d} 
		 && \text{on $ \partial \Solid^*_{\text{tail}}$},
	\end{aligned}
	\right.
\end{equation}
\normalsize
here $x^*_{cm}$ is the center of the swimmer's head $ \Solid^*_{\text{head}}$. 
To describe the deformation of the hyper-elastic flagella, we consider the Saint-Venant-Kirchhoff model. 
In this framework, the inertial and surface loads are defined in terms of the deformation 
gradient 
$$
F(\eta^t) = \mathbb{I}_d + \nabla \eta^t,
$$
and the second Piola-Kirchhoff tensor 
$$
\Sigma(\eta^t) = \lambda Tr(\epsilon(\eta^t)) \mathbb{I}_d + 2 \mu \epsilon(\eta^t),
$$
with $\lambda, \mu$ the Lamé coefficients and $\epsilon(\eta^t)$ the Green-Lagrange tensor. 

\subsection{Arbitrary Lagrangian-Eulerian framework} \label{model:ALE}

The Navier-Stokes equations are formulated in an Eulerian reference frame, where they describe the evolution of hydrodynamics within a fixed computational domain. However, the dynamics of the swimmers are described in a Lagrangian reference frame.  
To account for the geometric coupling condition, we adopt the Arbitrary Lagrangian-Eulerian framework \cite{chabannes_vers_2013}. 
This approach allows the fluid domain to follow the motion of the fluid-structure interaction interface. 
Specifically, the reference frame remains Lagrangian near the swimmers but stays Eulerian farther from the fluid-structure interaction interface.  
Both the swimmer and the fluid are discretized, and the discretized fluid domain follows the motion of the swimmer.

We define the ALE map $\Alemap^t_{\Fluid^*} : \, ]0,T] \times \Fluid^* \to \Fluid^t$ 
as a continuous and bijective function that gives the position of a particle in the current fluid domain $\Fluid^t$ based on its position in the reference domain $\Fluid^*$ 
$$
\Alemap^t_{\Fluid^*}(x^*) = x^t.
$$
To construct this ALE map, we introduce the displacement of the discretized fluid domain over time, 
denoted as $\eta_{\Fluid}^t :  \, ]0,T] \times \Fluid^* \to \Fluid^t$. A possible definition of the ALE map is then:  
$$
\Alemap^t_{\Fluid^*}(x^*) = x^* + \eta_{\Fluid}^t(x^*).
$$  
The fluid domain displacement $\eta_{\Fluid}^t$ can be obtained, for instance, by a harmonic extension
 of the swimmer’s displacement at time $t$, defined on $\partial \Solid^*$, to the interior 
of the fluid domain. This extension is determined by solving the following Laplace smoothing equation

\small
\begin{equation}
	\left\{
	\begin{aligned}
		\nabla \cdot\big((1+\tau) \nabla \eta_{\Fluid}^t \big) &= 0_{\Real^d} 
		\, && \text{in $ \Fluid^*$},\\
		\eta_{\Fluid}^t &= \eta^t 
		\, && \text{on $ \partial \Solid_i^*$},
	\end{aligned}
	\right.
    \label{Eq:ALE}
\end{equation}
\normalsize
\noindent where $\eta^t$ is the solution of \eqref{eq:hyper-elasticity}, and $\tau$ acts as a 
space-dependent diffusion coefficient. 
 It is a piecewise constant coefficient, defined on 
each element \( e \) of the domain's discretization as  
\[ \tau\big |_e = \frac{1 - V_{\min}/V_{\max}}{V_e/V_{\max}}, \]  
where \( V_{\max} \), \( V_{\min} \), and \( V_e \) are the volumes of the largest, smallest, and current elements of the domain discretization, 
respectively \cite{kanchi_3d_2007}. This coefficient \( \tau \) ensures that mesh deformation is applied primarily to elements with larger volumes.

When transitioning from the fluid Eulerian frame to an ALE frame, 
an additional term appears in the Navier-Stokes equations to account for the velocity 
of the moving domain, denoted by $u_{\mathcal{A}}^t$. 
The full system in the ALE reference frame is divided into two problems.
First, the \textit{fluid–rigid} problem, solving the hydrodynamics and the swimmer 
rigid motion, is governed by the following equations 
\small
\begin{equation}
	\left\{
	\begin{aligned}
		\rho_\Fluid \Big(\partial_t u^t  + \big( (u^t - u_{\mathcal{A}}^t) \cdot \nabla \big) u^t \Big) - \nabla \cdot \sigma(u^t,p^t) &= 0_{\Real^d}
		\, && \text{in $ \Fluid^t$},\\
		\nabla \cdot u^t &= 0 
		\, && \text{in $ \Fluid^t$},\\
		u^t &= 0_{\Real^d}
		\, && \text{on $ \partial \Fluid^t_D$},\\
		\sigma(u^t,p^t) n_{\Fluid}^t &=  0_{\Real^d} 
		\, && \text{on $ \partial \Fluid^t_N$},\\
		u^t &= \tilde{u}^t 
		\, && \text{on $ \partial \Solid^t$},\\
		m^t \, \mathrm{d}_{t} U^t &= f^t_H,\\
		\mathrm{d}_{t} \Big[R(\theta^t) J^t R(\theta^t)^T \omega^t\Big] &= T_H^t + T^t_m,
	\end{aligned}
	\right.
    \label{Eq:fluid_rigid:new}
\end{equation}
\normalsize
where $\tilde{u}^t =  U^t + \omega^t \times (x^t - x^t_{cm}) + \partial \eta^t  \circ (\mathcal{A}^t_{\Fluid^*})^{-1}$ 
the velocity of the swimmer.
Then, the \textit{fluid–elastic} problem computes the elastic deformations and 
introduces the continuity of displacements and stresses. This problem is described by
\small
\begin{equation}
	\left\{
	\begin{aligned}
		\rho_{\Solid} \, \partial_{tt} \eta^{t} - \nabla \cdot \big( F(\eta^t) \, \Sigma(\eta^t) \big) &= 0_{\Real^d}  
		\, && \text{in $ \Solid^*$},\\
		\eta^{t} &= \eta_{T}^t + \eta_{R}^t 
		\, && \text{on $ \Solid_{\text{head}}^*$},\\
        F(\eta^t) \, \Sigma(\eta^t) \, n^*_{\Solid} &= - \bigl( \sigma(u^t,p^t) \, n_{\Fluid}^t \bigr) \circ \mathcal{A}^t_{\Fluid^*}  \, && \text{on $ \partial \Solid^*$},
	\end{aligned}
	\right.
    \label{Eq:fluid_elastic:new}
\end{equation} 
\normalsize
where the translational and rotational displacements are respectively given by
$$
\mathrm{d}_{t} \eta_{T}^t = U^t, \quad \eta_{R}^t  = R(\theta^t)(x^* - x^*_{cm}) - (x^* - x^*_{cm}).
$$



\section{Numerical discretization}
\label{discretization}

For the spatial discretization of the discretized domain $\Omega_h^t = \Fluid_h^t \cup \Solid^t_h$, 
we use standard Lagrange finite elements. 
The associated discrete approximation spaces are 
constructed from piecewise polynomial functions of degree $N$, denoted $P_c^N(\Omega_h^t)$
\begin{align*}
P_c^N(\Omega_h^t) &= \left\{ v \in C^0(\Omega_h^t) \mid v \circ \varphi_{K_e}^{\mathrm{geo}} \in \mathbb{P}^N(\widehat{K}), \quad \forall K_e \right\}, 
\end{align*}
where $\mathbb{P}^N(\widehat{K})$ denotes the space of scalar polynomials of total
 degree less than or equal to $N$ defined on the reference element $\widehat{K}$, 
 and $\varphi_{K_e}^{\mathrm{geo}} \, : \widehat{K} \to K_e$ the geometric transformation 
 allowing to obtain each element $K_e$ from the reference element $\widehat{K}$. 
Vector-valued finite element spaces are constructed as Cartesian products of 
the scalar-valued spaces, denoted by
$
[P_c^N(\Omega_h^t)]^d.
$
We also use the notation $H^m(\Omega_h^t)$ to denote the Sobolev space of order 
$m$, and with boundary condition $g$ on a subset 
$\Gamma \subset \partial \Omega_h^t$ it is given by
$$
H^m_{(g,\Gamma)}(\Omega_h^t) = \{ v \in H^m(\Omega_h^t) \mid v|_{\Gamma} = g \}.
$$
The vector-valued versions are respectively given by $[H^m(\Omega_h^t)]^d$
and $[H^m_{(g,\Gamma)}(\Omega_h^t)]^d$.

Regarding temporal discretization, the time interval 
$]t_0 = 0, T = t_{N_t}]$, with $T > 0$ denoting the final time, 
is divided into $N_t$ time steps of size $\delta t$. The discrete time 
instances are denoted by $t_n$ for $n = 1, \dots, N_t$.
Finally, a mesh at time $t_{n+1}$ is denoted by $\Omega_h^{n+1}$, 
and a discrete solution at this time is written as $u_h^{n+1}$. 

We will first present the spatial and temporal discretization of the ALE map 
before introducing that of the \textit{fluid–rigid} and \textit{fluid–elastic} problems.

\subsection{Discretization of the ALE map}

The approximation spaces at time $t_{n+1}$ of the trial and test functions are respectively given by: 
\small
\begin{align*}
	V_{\mathcal{A},h}^{*} &= \{  \mathbf{v}, \mathbf{v}  \in [H^1_{(0_{\Real^d}, \partial \Fluid^*_h)}(\Fluid^*_h)]^d \cap [P^1_c(\Fluid^*_h)]^d , \mathbf{v} =  \eta_{\Solid^*,h}^{n+1} \text{ on } \partial \mathcal{S}_h^* \},\\
	W_{\mathcal{A},h}^{*} &= \{  \mathbf{v}, \mathbf{v}  \in [H^1_{(0_{\Real^d}, \partial \Fluid^*_h)}(\Fluid^*_h)]^d \cap [P^1_c(\Fluid^*_h)]^d , \mathbf{v} =  0_{\Real^d} \text{ on } \partial \mathcal{S}_h^* \}.
\end{align*}
\normalsize
Thus, the fluid domain displacement is discretized using continuous affine finite elements. 
The space- and time-discretized version of \eqref{Eq:ALE} at time $t_{n+1}$ gets

\medskip 

Find $\eta_{\Fluid,h}^{n+1} \in V_{\mathcal{A},h}^{*}$ such that
\begin{equation}\label{eq:ALE_EDP_dis}
\underbrace{\int_{\Fluid_h^*} (1+\tau)\nabla \eta_{\Fluid,h}^{n+1} : \nabla \mathbf{v}}_{:= a_{\mathcal{A}^*}(\eta_{\Fluid^*,h}^{n+1},\mathbf{v})} = 0 \quad \text{for all } \mathbf{v} \in W_{\mathcal{A},h}^{*}.
\end{equation}

\noindent The discrete ALE map is then defined as  
\[
\Alemap^{n+1}_{\Fluid^*,h}(x^*) = x^* + \eta_{\Fluid,h}^{n+1}(x^*).
\]
Large deformations may lead to element inversions, resulting in invalid triangulation. 
To avoid this, mesh quality metrics \cite{field_qualitative_2000} are used. 
If the quality remains above a fixed threshold, the domain is updated using the 
ALE map. Otherwise, remeshing is performed before applying it. 
Details of these procedures are given in \cite{van_landeghem_towards_2024}.

For the numerical resolution, we define the algebraic system associated with \eqref{eq:ALE_EDP_dis} as follows
\[
A \boldsymbol{\eta_{\Fluid}} = 0_{\Real^{N_{\boldsymbol{\eta_{\Fluid}}}}},
\]
where we denote by
\[
A = \left(a_{\mathcal{A}^*}(\xi_i, \xi_j)\right)_{i,j} \in \Real^{N_{\boldsymbol{\eta_{\Fluid}}} \times N_{\boldsymbol{\eta_{\Fluid}}}}, \quad \text{for } 1 \leq i,j \leq N_{\boldsymbol{\eta_{\Fluid}}},
\]
and $(\xi_i)_i$, $i = 1, \cdots, N_{\boldsymbol{\eta_{\Fluid}}}$  is the basis of the space $V_{\mathcal{A},h}^{*}$.
The solver and preconditioner applied to this algebraic system are detailed in \ref{strategies}.

\subsection{Discretization of the fluid-rigid problem}
\label{framework::fluid}
The Dirichlet boundary condition, enforcing the coupling of velocities at the 
interface of the swimmer is incorporated both in the function space associated with 
the fluid velocity and in the space of test functions.  
The discrete spaces for the fluid velocity and pressure at time~$t_{n+1}$ are then respectively defined as
\small
\begin{align*}
V^{n+1}_{\Fluid,h} 
  &= \Bigl\{ 
      \mathbf{v} : \Fluid^{n+1}_h \to \Real^d,
      \mathbf{v} = \mathbf{\tilde{v}} \circ (\mathcal{A}^{n+1}_{\Fluid^*,h})^{-1},\ 
      \mathbf{\tilde{v}} \in [H^1_{(0_{\Real^d},\mathcal{F}^*_{D,h})}(\Fluid^*_h)]^d 
      \cap [P_c^M (\Fluid^{*}_h)]^d, \\
  &\quad \quad 
      \mathbf{\tilde{v}} = \mathbf{\tilde{u}} \text{ on } \partial \Solid^{n+1}_h,\ 
      \mathbf{\tilde{u}} \in \Real^d 
    \Bigr\}, \\[1mm]
Q^{n+1}_{\Fluid,h} 
  &= \Bigl\{ 
      \mathbf{q} : \Fluid^{n+1}_h \to \Real,\ 
      \mathbf{q} = \mathbf{\tilde{q}} \circ (\mathcal{A}^{n+1}_{\Fluid^*,h})^{-1},\
      \mathbf{\tilde{q}} \in P_c^N (\Fluid^{*}_h)
    \Bigr\}.
\end{align*}
\normalsize
We choose the Taylor-Hood inf-sup finite elements to satisfy the compatibility condition. 
The polynomial degrees are related by $N = M - 1$, and typically we set $M = 2$ and 
$N = 1$. Thus, the velocity and pressure are respectively discretized using continuous 
piecewise quadratic, and affine finite elements.  
Following the ideas of \cite{maury_direct_1999}, the discretized formulation of the \textit{fluid–rigid} 
problem~\eqref{Eq:fluid_rigid:new} deriving from the variational form is stated as follows

\medskip

Find the fluid velocity and pressure 
\( (u^{n+1}_{h}, p^{n+1}_{h}) \in V^{n+1}_{\Fluid,h} \times Q^{n+1}_{\Fluid,h} \), and the swimmer's linear and angular velocities $U^{n+1}, \omega^{n+1} \in \Real^d \times \Real^{d^*}$, such that for all \( (\mathbf{v}, \mathbf{q}) \in V^{n+1}_{\Fluid,h} \times Q^{n+1}_{\Fluid,h} \), 
and $\mathbf{U}, \boldsymbol{\omega} \in \Real^d \times \Real^{d^*} $, one has 

\small 
\begin{equation}
\begin{aligned}
a_{\Fluid^t}\bigl(u^{n+1}_{h},\mathbf{v}\bigr)
  - b_{\Fluid^t}\bigl(\mathbf{v},p^{n+1}_{h}\bigr) 
+ m^{n+1}\,\mathrm{d}_{t}U^{n+1} \cdot \mathbf{U} \\ 
\qquad
 + \, \mathrm{d}_{t}\Bigl[\,R(\theta^{n+1})\,J^{n+1}\,R(\theta^{n+1})^{T}\,\omega^{n+1}\Bigr] \cdot \boldsymbol{\omega} 
&= l_{\Fluid}\bigl(\mathbf{v},\mathbf{U},\boldsymbol{\omega}\bigr),\\[2mm]
b_{\Fluid^t}\bigl(u^{n+1}_{h},\mathbf{q}\bigr) & = 0.
\end{aligned}
\label{Eq:fluid_rigid}
\end{equation}
\normalsize
with the bilinear and linear forms
\small
\begin{align*}
	a_{\Fluid^t} (u^{n+1}_{h}, \mathbf{v}) 
	&= \rho_{\Fluid}  \int_{\Fluid^{n+1}_h} \partial_t u_{h}^{n+1}  \cdot \mathbf{v}  
	+ \rho_{\Fluid} \int_{\Fluid^{n+1}_h} \left( (u^{n+1}_{h} -  u_{\mathcal{A},h}^{n+1})
	 \cdot \nabla \right) u^{n+1}_{h} \cdot \mathbf{v} \\
	&\quad + \mu_{\Fluid} \int_{\Fluid^{n+1}_h} \nabla u_{h}^{n+1} :  \nabla \mathbf{v}, \\
	l_{\Fluid^t} (\mathbf{v}) 
	&= f_H^{n+1} \cdot \mathbf{U} + T_H^{n+1} \cdot \boldsymbol{\omega} + T_m^{n+1} \cdot \boldsymbol{\omega}  , \\
	b_{\Fluid^t}(u^{n+1}_{h}, \mathbf{q}) 
	&= - \int_{\Fluid^{n+1}_h} \mathbf{q} \cdot \nabla \cdot u^{n+1}_{h}.
\end{align*}
\normalsize
The time derivative of $u_h^{n+1}$ as well as the ALE velocity $u^{n+1}_{\mathcal{A},h}$ 
are time-discretized using the second-order backward differentiation formula, BDF2. 
For the numerical resolution, the degrees of freedom of the fluid velocity on 
the interface of the swimmer are treated differently since they depend only on the linear 
and angular velocities. 

We denote the vectors 
$\mathbf{u} = (\mathbf{u}_1, \dots, \mathbf{u}_{N_{\mathbf{u}}}) = (\mathbf{u_I}, \mathbf{u_\Gamma})$ for the velocity, 
where the subscript $\Gamma$ represent the degrees of freedom associated with the interface,  
$\mathbf{p} = (\mathbf{p}_1, \dots, \mathbf{p}_{N_{\mathbf{p}}})$ for the pressure, and 
$\mathbf{U} \in \mathbb{R}^d$ and $\boldsymbol{\omega} \in \mathbb{R}^{d^*}$ 
represent the linear and angular velocities, respectively. 
The resulting algebraic system of \eqref{Eq:fluid_rigid} is then given by
\begin{equation}
    \begin{bmatrix}
        A_{II} & A_{I\Gamma} & 0 & 0 & B_{I}^T  \\
        A_{\Gamma I} & A_{\Gamma \Gamma}  & 0 & 0 & B_{\Gamma}^T\\
        0 & 0 & m^{n+1} \mathbb{I}_d
        & 0 & 0 \\
        0 & 0 & 0& R(\theta^{n+1}) J^{n+1} R(\theta^{n+1})^T
        & 0\\
        B_I & B_\Gamma & 0 & 0 & 0   \\
    \end{bmatrix}
    \begin{bmatrix}
        \mathbf{u_I} \\ \mathbf{u_{\Gamma}} \\ \mathbf{U} \\ \boldsymbol{\omega} \\ \mathbf{p}
    \end{bmatrix} 
    = 
    \begin{bmatrix}
        0 \\
        0\\
        l_U
        \\
        l_{\omega}
        \\
        0
    \end{bmatrix},
    \label{Eq:flu-rig-block}
\end{equation} 
where
\begin{align*}
A_{JK} &= (a_{\Fluid^t}(\phi_{J_{i}},\phi_{K_{j}}))_{i,j} \in \Real^{N_{\mathbf{u_J}} \times N_{\mathbf{u_K}}},   \textrm{ for } J,K \in \{I,\Gamma\},\\
B_J &= (b_{\Fluid^t}(\phi_{J_{i}},\psi_{j}))_{i,j}  \in \Real^{N_{\mathbf{u_J}} \times N_{\mathbf{p}}} \textrm{ for } J \in \{I,\Gamma\} ,\\
l_{U} &= f_H^{n+1} \cdot \mathbf{U} , \\
l_{\omega} &= T_H^{n+1} \cdot \boldsymbol{\omega} + T_m^{n+1} \cdot \boldsymbol{\omega},
\end{align*}
and $(\phi_i)_i$, $i = 1, \cdots, N_{\mathbf{u}}$, respectively $(\psi_i)_i$, $i = 1, \cdots, N_{\mathbf{p}}$, 
the basis of the spaces $V^{n+1}_{\Fluid,h}$ and $Q^{n+1}_{\Fluid,h}$. 

Using the no-slip boundary conditions
$\mathbf{u} = \mathbf{\tilde{u}} = \mathbf{U} + \boldsymbol{\omega} \times (x^{n+1}-x^{n+1}_{cm})$ on  $\partial \Solid^{n+1}_h$, 
we introduce the operator $\mathcal{P}$ such as 
\begin{equation}
\label{eq:boundary_constraint}
 (\mathbf{u_I},\mathbf{u_{\Gamma}},\mathbf{U},\boldsymbol{\omega},\mathbf{p} )^T = \mathcal{P} (\mathbf{u_I},\mathbf{U},\boldsymbol{\omega},\mathbf{p} ),
\end{equation}
with 
\begin{equation*}
    \mathcal{P} =
    \begin{bmatrix}
        \mathbb{I}_d &  0 & 0 & 0\\
        0 &\tilde{P}_{U} & \tilde{P}_{\omega} & 0\\
        0 &\mathbb{I}_d & 0 & 0\\
        0 &0 & \mathbb{I}_d & 0\\
        0 &0& 0& \mathbb{I}_d
    \end{bmatrix}\,.
\end{equation*}

Finally, by using \eqref{eq:boundary_constraint}, and by multiplying by $\mathcal{P}$, 
one obtains the algebraic system of the fluid-swimmer problem
\small 
\begin{equation*}
    \mathcal{P}^T \begin{bmatrix}

            A_{II} & A_{I\Gamma} & 0 & 0 & B_{I}^T  \\
            A_{\Gamma I} & A_{\Gamma \Gamma}  & 0 & 0 & B_{\Gamma}^T\\
            0 & 0 & m^{n+1} \mathbb{I}_d
            & 0 & 0 \\
            0 & 0 & 0& R(\theta^{n+1}) J^{n+1} R(\theta^{n+1})^T
            & 0\\
            B_I & B_\Gamma & 0 & 0 & 0   \\

    \end{bmatrix} \mathcal{P} \begin{bmatrix}
        \mathbf{u_I} \\ \mathbf{U} \\ \boldsymbol{\omega} \\ \mathbf{p}
    \end{bmatrix}  
 = 
 \mathcal{P}^T 	
 \begin{bmatrix}
    0 \\
    0\\
    l_{U}
    \\
    l_{\omega}
    \\
    0
    \end{bmatrix}. 
    \label{algebraicsystem::fluid:rigid}
\end{equation*}
\normalsize

\subsection{Discretization of the fluid–elastic problem}

The approximation spaces of admissible elastic displacements, and of the tests 
functions, at time $t_{n+1}$, are respectively defined by
\small
\begin{align*}
	V_{\Solid,h}^{*} &= \{  \mathbf{v}, \mathbf{v}  \in [H^1_{(\eta_T^{n+1} + \eta_R^{n+1}, \partial \Solid^*_{\text{head},h})}(\Solid^*_h)]^d \cap [P^1_c(\Solid^*_h)]^d  \},\\
	W_{\Solid,h}^{*} &= \{  \mathbf{v}, \mathbf{v}  \in [H^1_{(0_{\Real^d}, \partial \Solid^*_{\text{head},h})}(\Solid^*_h)]^d \cap [P^1_c(\Solid^*_h)]^d \}.
\end{align*}
\normalsize
The fully discretized \textit{fluid–elastic} problem is stated as follows

\medskip 

Find the elastic displacement $\eta^{n+1}_{h} \in V^{*}_{\Solid,h}$ such that, for all 
$\mathbf{v} \in W^{*}_{\Solid,h}$, one has
\small
\begin{align*}
    a_{\Solid^*}(\eta_{h}^{n+1},\mathbf{v}) 
    = - \int_{\partial \Solid_{h}^{n+1}} 
    |\det F(\eta_{h}^{n+1})| 
    \, ||(F(\eta_{h}^{n+1}))^{-T}|| 
    \bigl( \sigma (u_h^{n+1}, p_h^{n+1}) n_{\Fluid}^{n+1} \bigr) \cdot \mathbf{v},
\end{align*}
\normalsize
where the boundary condition is expressed on the current interface 
$\partial \Solid_{h}^{n+1}$, since the hydrodynamic forces and torques are resolved 
onto the current domain, see \eqref{eq:hydro_forces}.
The bilinear form is defined as
\begin{align*}
    a_{\Solid^*}(\eta_{h}^{n+1},\mathbf{v}) 
    &= \rho_{\Solid}
    \int_{\Solid^*_h} \partial_t \eta^{n+1}_{h} \cdot \mathbf{v} 
    + \int_{\Solid^*_h} F(\eta_h^{n+1}) \Sigma(\eta_h^{n+1}) : \nabla \mathbf{v}.
\end{align*}
Once again, we consider the BDF2 scheme to discretize the time derivative. 
The corresponding algebraic system is written as
\begin{align}
A \boldsymbol{\eta_{\Solid}} = G,
\label{eq:algebraic:elasticity}
\end{align}
where the stiffness matrix $A \in \mathbb{R}^{N_{\boldsymbol{\eta_{\Solid}}} \times N_{\boldsymbol{\eta_{\Solid}}}}$ 
is expressed in the basis $(\varphi_i)_i$, $i = 1, \cdots, N_{\boldsymbol{\eta_{\Solid}}}$ 
of the discrete space $V^{*}_{\Solid,h}$ as
\[
A = \bigl(a_{\Solid^*}(\varphi_i, \varphi_j)\bigr)_{i,j}, 
\]
and the right-hand side vector $G \in \mathbb{R}^{N_{\boldsymbol{\eta_{\Solid}}}}$ 
includes the boundary conditions.



\section{Resolution strategies}
\label{strategies}

\subsection{Full algorithm} 

The complete coupled algorithm for the elasto-magneto-swimmer is given in \cref{algo}.
At each iteration, the swimmer's displacement is updated using a relaxation method, where the relaxation parameter $t^{k+1}$ is computed using the Aitken method, as described in \cite{kuttler_fixed-point_2008}.

\medskip 

\begin{algorithm}[!ht]
\caption{Fixed-Point algorithm for the magneto-swimmer}
\begin{small} 
\begin{algorithmic}[1]
\STATE \textbf{Input:} Solid displacements at previous time steps $\eta_h^{n-1}$ and $\eta_h^n$, tolerance $\mathrm{tol}$, maximum iterations $k_{\text{max}}$.
\STATE \textbf{Output:} New solid displacement $\eta_h^{n+1}$.
\STATE Initialize iteration counter: $k = 0$.
\STATE Predict initial displacement using BDF2-consistent extrapolation 
$$\eta_h^{n+1,0} = 2 \eta_h^n - \eta_h^{n-1}.$$
\While{$\epsilon > \mathrm{tol}$ and $k < k_{\text{max}}$}
    \STATE Compute ALE map $\tilde{\ALE}_{\Fluid^*,h}^{n+1,k+1}$ from $\eta_h^{n+1,k}$.\label{step1}
    \STATE Solve \textit{fluid–rigid} problem on $\tilde{\ALE}_{\Fluid^*,h}^{n+1,k+1}(\Fluid^*_h)$ to obtain $\eta_T^{n+1,k+1} + \eta_R^{n+1,k+1}$. \label{step2}
\STATE Update total solid displacement on the reference interface 
$$\tilde{\eta}_h^{n+1,k+1} = \eta_h^{n+1,k} + \eta_T^{n+1,k+1} + \eta_R^{n+1,k+1}.$$
\STATE Compute ALE map $\ALE_{\Fluid^*,h}^{n+1,k+1}$ from $\tilde{\eta}_h^{n+1,k+1}$.
    \STATE Solve \textit{fluid–elastic} problem on $\Solid^*_h$ to obtain total displacement $\tilde{\tilde{\eta}}_h^{n+1,k+1}$.\label{step4}
    \STATE Compute error $$\epsilon = \|\tilde{\tilde{\eta}}_h^{n+1,k+1} - \eta_h^{n+1,k}\|_{L^2(\Solid^*_h)}.$$
    \IF{$\epsilon > \mathrm{tol}$}
        \STATE Update swimmer displacement $$\eta_h^{n+1,k+1} = t^{k+1} \tilde{\tilde{\eta}}_h^{n+1,k+1} + (1-t^{k+1}) \eta_h^{n+1,k}.$$
        \STATE Subtract rigid motion $$\eta_h^{n+1,k+1} = \eta_h^{n+1,k+1} - (\eta_T^{n+1,k+1} + \eta_R^{n+1,k+1}).$$
        \STATE Increment iteration counter $k = k+1$.
    \ELSE
        \STATE Set $\eta_h^{n+1} = \tilde{\tilde{\eta}}_h^{n+1,k+1}$ and exit loop.
    \ENDIF
\EndWhile
\end{algorithmic}
\end{small}
\label{algo}
\end{algorithm}

\subsection{Remeshing strategy}

The mesh deformation described by the ALE map may lead to significant element 
distortion during long-term simulations. 
To preserve mesh quality, a remeshing procedure is performed at a prescribed 
frequency, treated as a user-defined parameter depending on the simulated trajectories. 
At each remeshing step, the reference domain is updated accordingly, and all 
fields are projected onto the new mesh. 
The remeshing operations are carried out using the \texttt{MMG} tool in 
sequential mode, or \texttt{ParMMG} in parallel \cite{balarac_tetrahedral_2022}.
The latter libraries take as input the current domain, which could be of bad quality, 
and a scalar function defined as the remeshing metric.
This metric specifies the desired characteristic size of the elements in the resulting mesh.
To accurately capture the fluid-swimmer interaction, we adopt a graded remeshing strategy, 
meaning that the mesh size is adjusted according to the distance to the 
swimmer’s boundary in the current fluid domain.

\subsection{Solver and preconditioner}
\label{prec::elasticity}

The numerical scheme of the magneto-swimmer model is based on solving a 
fixed-point algorithm which iterates between the ALE \eqref{Eq:ALE}, \textit{fluid-rigid} \eqref{Eq:fluid_rigid:new}, and 
\textit{fluid-elastic} \eqref{Eq:fluid_elastic:new} problems corresponding to the step \ref{step1}, \ref{step2} and \ref{step4} of the full algorithm.
Each sub-problem is solved using appropriate solvers and preconditioners, 
which are summarized for both two- and 
three-dimensional cases in \cref{tab:preconditionner:magneto}. 
The $\feelpp$ library interfaces with PETSc \cite{argonne_national_laboratory_petsc_nodate} 
for the efficient solution of large-scale linear and nonlinear systems. 

\small 
\begin{table}[!ht]
    \centering
    \begin{tabular}{lcc}
    \toprule
    \textbf{Sub-problem} & \textbf{Two Dimensions} & \textbf{Three Dimensions} \\
    \midrule
    ALE map & CG + GAMG & CG + GAMG \\
    \textit{Fluid-rigid} & Direct solver + LU & GMRES + Block preconditioner\\
    \textit{Fluid-elastic} & Direct solver + LU & GMRES + GASM \\
    \bottomrule
    \end{tabular}
    \caption{Solvers and preconditioners for the magneto-swimmer model.}
    \label{tab:preconditionner:magneto}
\end{table}
\normalsize



\section{Numerical results}
\label{results}

We perform different numerical tests to validate and illustrate the proposed magneto-elastic swimmer model.
First, we analyze the relationship between the swimmer's net displacement and the area enclosed by its stroke in the configuration space.
Second, we investigate the net displacement of the swimmer as a function of the external magnetic field frequency in two dimensions.
Finally, we extend the validation to three dimensions by simulating the full 3D magneto-swimmer and comparing the results with the 2D case and available literature. 
The main physical parameters used in the simulations are summarized in  
\cref{tab::properties_swimmer}, these values are chosen according to \cite{oulmas_comparing_2019}. We prescribe the external magnetic field and the swimmer's magnetic moment based on experimental data from \cite{oulmas_3d_2017,oulmas_comparing_2019},
\begin{equation}
B^{n} = (b_x, b_y \sin(2\pi f t_{n}))^T,
\label{eq:magnetic_field}
\end{equation}
with $b_x = b_y = 0.005$T. 
Under the action of this oscillating magnetic field, the swimmer follows a straight horizontal path.

\begin{table}[!ht]
	
	\begin{center}
		\begin{tabularx}{0.9\textwidth} { 
			| >{\raggedright\arraybackslash}X 
			| >{\centering\arraybackslash}X | }
		    \hline  
		    Physical parameter & Value \\  
		    \hline 
		    Young modulus of the head &  $41$GPa \\  
		    \hline
		    Poisson coefficient of the head & $0.281$  \\ 
		    \hline 
            Magnetization of the head & $10^5 \text{A}/\text{m}$  \\ 
		    \hline 
		    Density of the head & $7000 \text{kg}/\text{m}^3$  \\ 
		    \hline 
            Height of the head & $0.5 \text{mm}$  \\ 
		    \hline
            Young modulus of the tail & $0.1$MPa   \\ 
			\hline 
			Poisson coefficient of the tail & $0.4$ \\ 
			\hline
            Density of the tail & $1300 \text{kg}/\text{m}^3$  \\ 
		    \hline 
            Length of the tail & $7.5 \text{mm}$  \\ 
		    \hline
            Maximal diameter of the tail & $1.5 \text{mm}$  \\ 
		    \hline
            Magnetic field intensity $b_x,b_y $ & $5  \text{mT}$  \\ 
		    \hline
	        \end{tabularx}
	
    \end{center}

	\caption{Physical properties of the magneto-swimmer taken from references 
    \cite{oulmas_3d_2017,oulmas_comparing_2019}.
    }

	\label{tab::properties_swimmer}
\end{table} 

\begin{figure}[!ht]
        \centering
        \begin{tikzpicture}[scale=1.75]

        \begin{scope}[shift={(0,0)}, rotate=10]

          \draw[black, thick] (0,0) -- (0.5,0);
          \draw[black, thick] (0,0) -- (0,1);
          \draw[black, thick] (0,1) -- (0.5,1);

          \draw[black, thick] (0.5,1) .. controls (1.8,0.4) .. (3,0.2);
          \draw[black, thick] (0.5,0) .. controls (1.8,-0.1) .. (3,0.2);

          \draw[->, blue, thick] (0.25,0.5) -- ++(-1,0);

          \draw[->, red, thick] (2.9,0.195) -- (3.5,0.24);

        \end{scope}

        \draw[dashed, thick, gray] (-0.1,0.53) -- ++(1.4,0);

        \draw[thin] (0.25,0.5) ++(0:0.5) arc[start angle=10, end angle=175, radius=0.5];
        \node at (-0.4,0.75) {\small{$\theta^n_{\text{head}}$}};

        \draw[thick, gray, dashed] (2.9,0.7) -- ++(0.8,0);

        \draw[thin] (2.9,0.7) ++(0:0.3) arc[start angle=0, end angle=18.3, radius=0.3];
        \node at (3.25,0.55) {\small{$\theta^n_{\text{tail}}$}};

        \end{tikzpicture}
    \caption{The illustration of the two angles used for the phase configuration space.}
    \label{fig:combined_swimmer}
\end{figure}
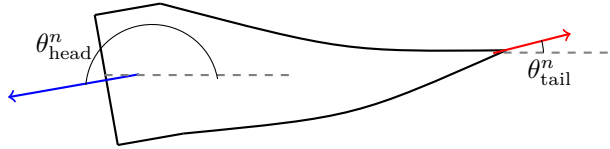

\begin{figure}[!ht]
    \centering
    \includegraphics[scale = 0.4]{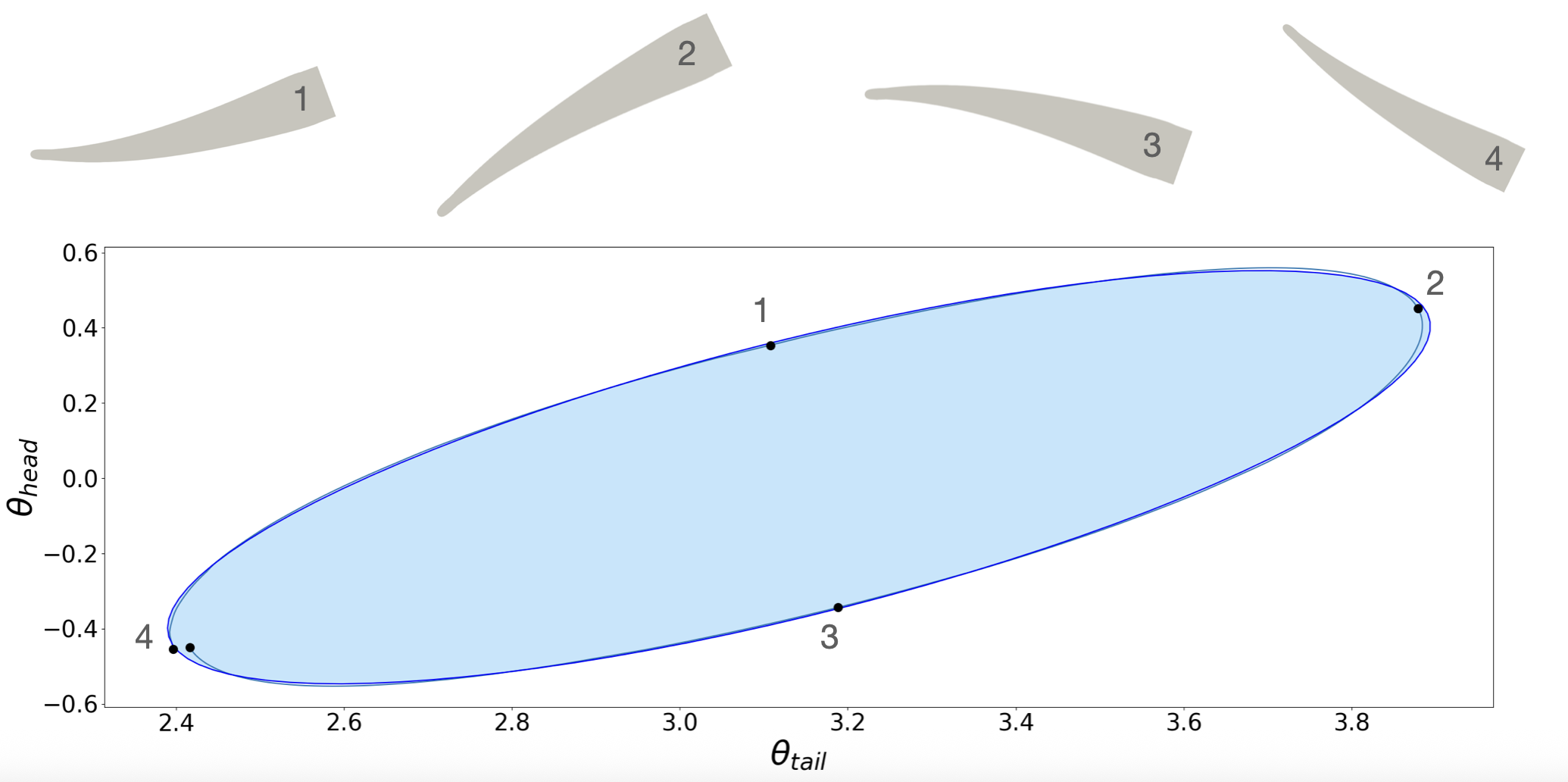}
\caption{At the top, the deformations of the swimmer with respect to the stroke cycle plotted at the bottom in the configuration space.}
 \label{fig:cycle_stroke}
\end{figure}

\paragraph{Net displacement vs swimmer's strokes}
We introduce the magneto-head and magneto-tail angles $\theta^n_{\text{head}}$ and $\theta^n_{\text{tail}}$ relative to the horizontal axis as illustrated in \cref{fig:combined_swimmer}.
We build their phase portrait based on their time evolution over one period as depicted in \cref{fig:cycle_stroke} to represent the swimmer's stroke in the configuration space.
The obtained curve approximates an ellipse, and according
to \cite{alouges_can_2015}, when the deformation is small, the mean $x$-displacement per period is directly correlated to the area of this ellipse. This result confirms the theoretical prediction that the swimmer will achieve an $x$-net displacement after the stroke cycle.

\paragraph{Net displacement vs frequency}

\begin{figure}[!ht]
	
	\centering  
	\includegraphics[width=\linewidth]{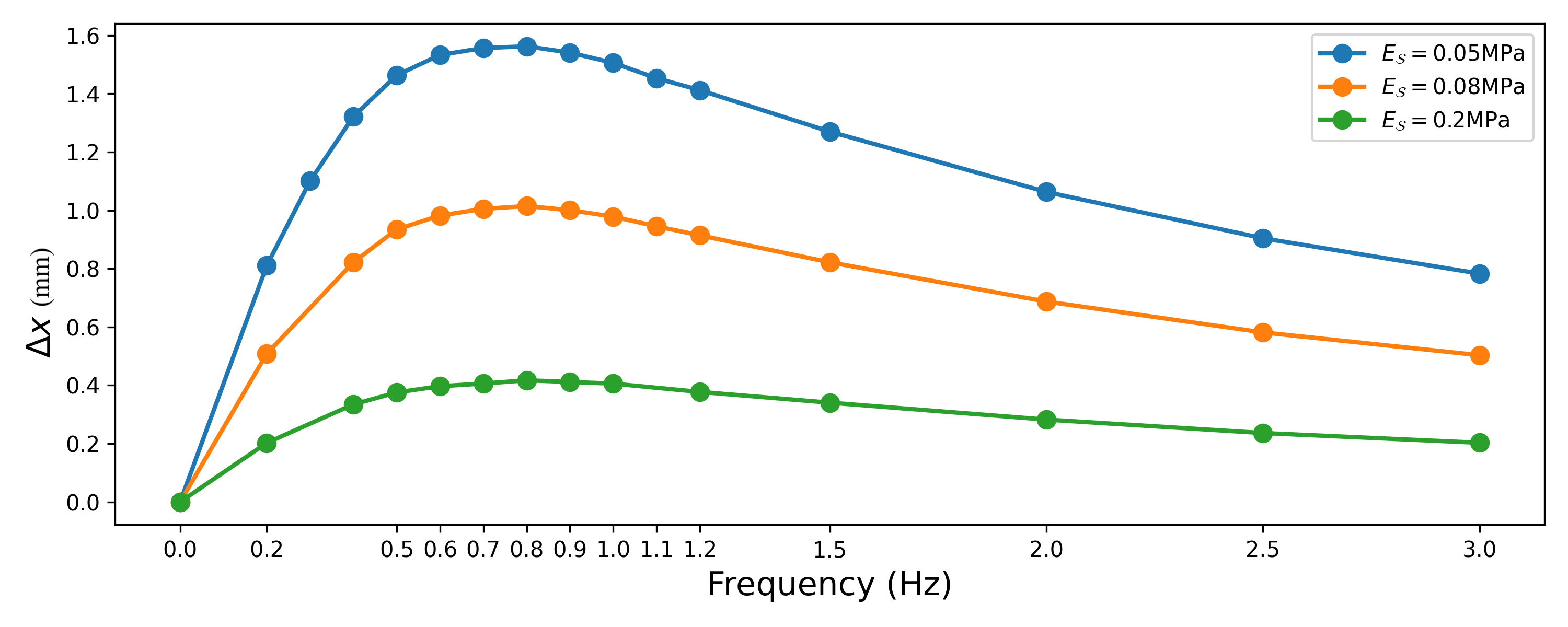}  
	
	\caption{Displacement of the magneto-swimmer over one period as a function 
	of the frequency of the magnetic field.}

	\label{Fig:freq_disp}
\end{figure}

In this case, we work with the two-dimensional model of the magneto-elastic swimmer.
In \cref{Fig:freq_disp}, we vary the frequency $f$ of the magnetic field defined in equation~\eqref{eq:magnetic_field} within the range 
$f \in [0\,\text{Hz}, 3\,\text{Hz}]$, while measuring the net displacement $\Delta x$ of the swimmer over one oscillation period of duration $\frac{1}{f}$. 
This analysis is performed for three different values of the tail’s Young modulus: 
$E_{\Solid} = 5 \times 10^4 \, \text{Pa}$, 
$E_{\Solid} = 8 \times 10^4 \, \text{Pa}$, 
and $E_{\Solid} = 2 \times 10^5 \, \text{Pa}$. 
The corresponding results are shown in \cref{Fig:freq_disp}, which plots the net displacement obtained during the third oscillation cycle, 
during the time interval $[\frac{2}{f}\, \text{s}, \frac{3}{f}\, \text{s}]$.  

We observe that the displacement per cycle increases with frequency, reaching a maximum around $f \approx 0.8\,\text{Hz}$, before decreasing as the frequency continues to rise. 
This non-monotonic trend is consistent with experimental and theoretical findings reported in the literature~\cite{alouges_can_2015,oulmas_3d_2017}, 
where an optimal frequency was also identified, corresponding to a resonance between the magnetic actuation and the elastic response of the swimmer’s tail. 
At low frequencies, the magnetic actuation is too slow to generate significant tail deformation, resulting in limited propulsion. 
Conversely, at high frequencies, the head cannot follow the rapid field oscillations, leading to attenuated elastic waves along the tail, reducing the net displacement.

Furthermore, we note that a smaller Young modulus leads to larger net displacements per cycle. 
This is attributed to the increased flexibility of the tail, which allows for greater deformation and thus more efficient conversion of oscillatory motion into forward propulsion. 
This observation also agrees with the results of~\cite{oulmas_3d_2017}, where flexible magnetic swimmers were shown to outperform stiffer ones under similar low-Reynolds-number conditions.

\paragraph{$3$D displacement of the elasto-magneto-swimmer}

\begin{figure}[!ht]
	
	\centering  
	\includegraphics[scale=0.5]{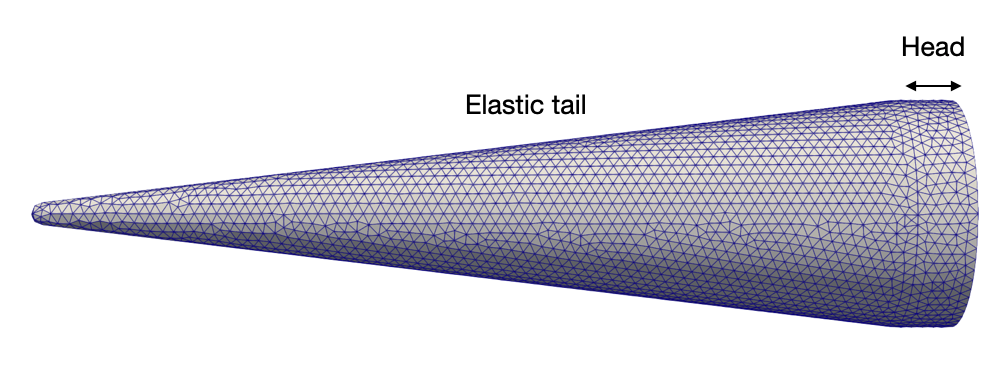}  
	
	\caption{Discretized geometry of the three-dimensional magneto-swimmer.}

	\label{Fig:geo_ms_3d_meshed}
\end{figure}

We now consider the same test case using the three-dimensional model of the elasto-magneto-swimmer. 
The geometry of the swimmer is shown in \cref{Fig:geo_ms_3d_meshed}. 
Its head is cylindrical, with a diameter of $1.5\,\text{mm}$ and a height of $0.5\,\text{mm}$, and is attached to a flexible tail of length $7.5\,\text{mm}$ and a minimum diameter of $0.2\,\text{mm}$. 
The density and Poisson’s ratio of both the head and the tail are set to typical values used in the literature for soft magnetic swimmers.\newline 

The Young’s modulus is fixed to $E_{\Solid} = 2 \times 10^5 \, \text{Pa}$, while 
the frequency of the external magnetic field is varied within the range 
$[0\,\text{Hz}, 3\,\text{Hz}]$. The external magnetic field in three 
dimensions is generalized by applying the same oscillatory pattern in the 
$y$-direction but keeping the $z$-component zero, see equations \eqref{eq:magnetic_field}. 
The intensity of the magnetic field is maintained at $b_x = b_y = 5\,\text{mT}$. 
The magnetization of the magneto-head remains fixed at $\mathsf{m} = 10^5 \text{A}/\text{m}$ 
as in the two-dimensional case and following experimental data from \cite{oulmas_3d_2017,oulmas_comparing_2019}.

\begin{figure}[!ht]
    \centering

    \includegraphics[width=1\textwidth]{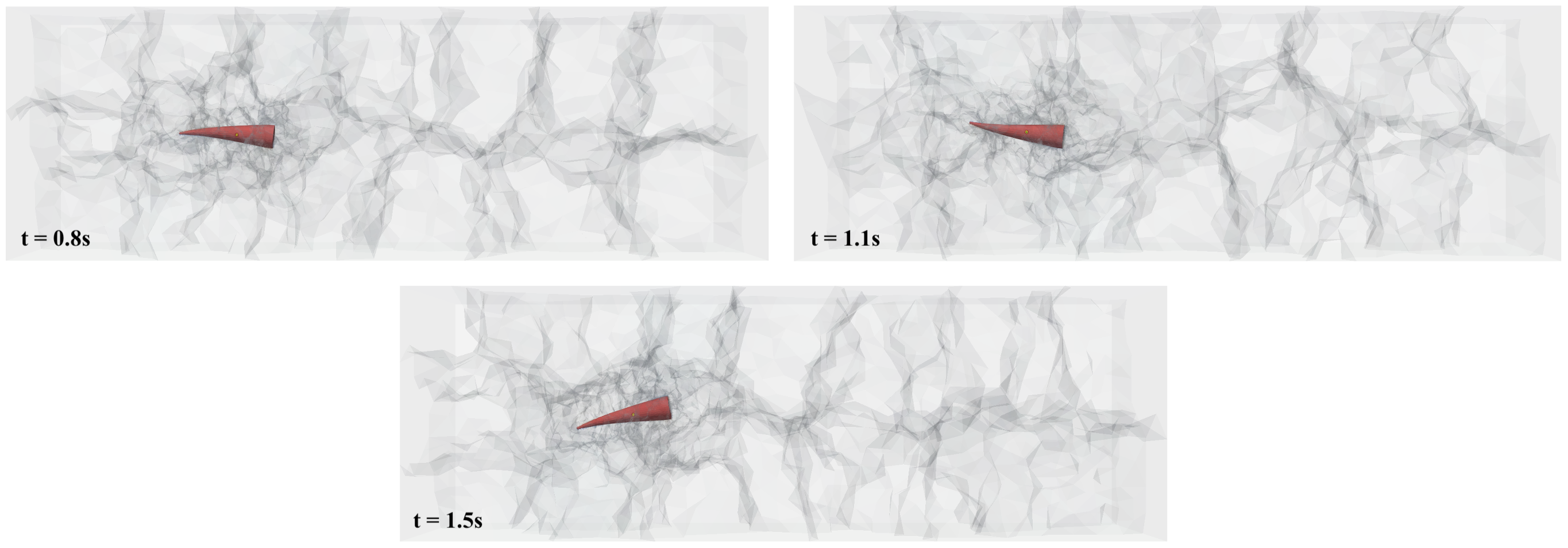}
    \caption{Visualization of the three-dimensional magneto-swimmer at different time instants. 
    Transparent gray surfaces represent the boundaries between domains used by the processors.}
    \label{fig:intro:magneto}
\end{figure}

\begin{figure}[!ht]
	\centering  
	\includegraphics[width=\linewidth]{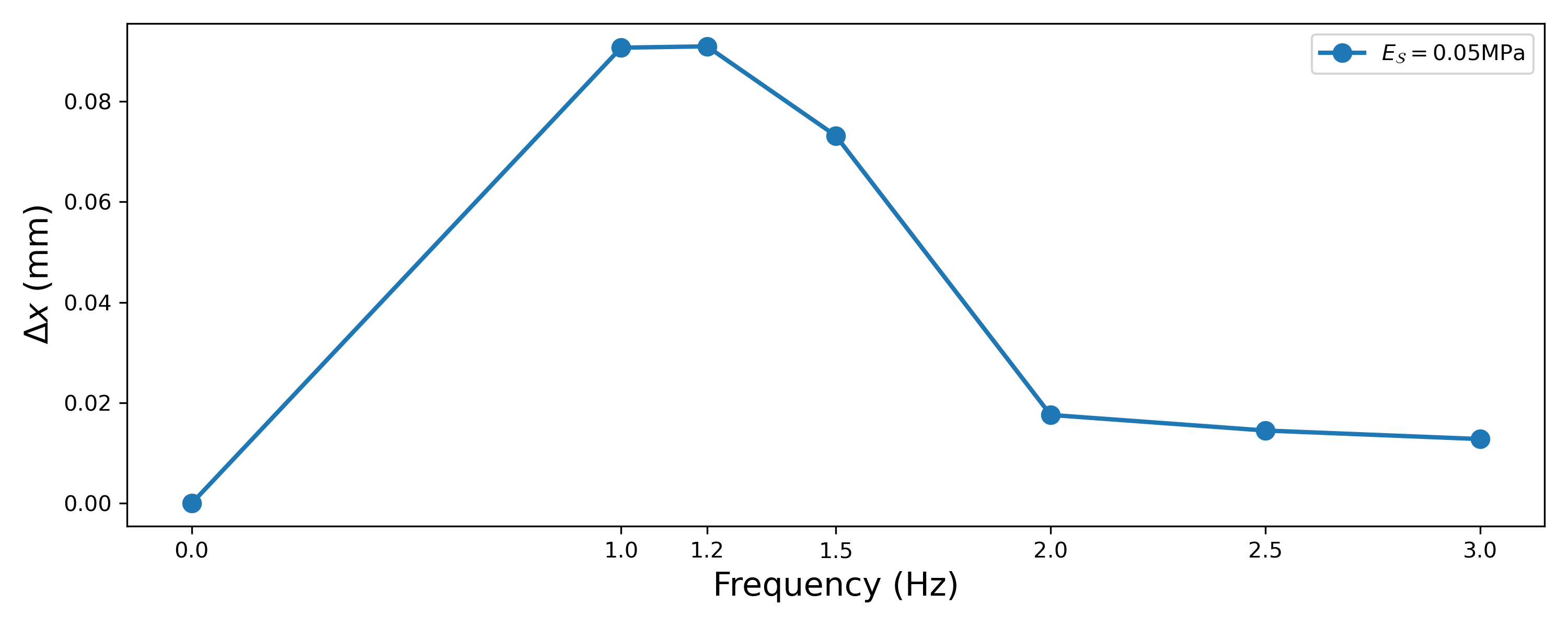}  
	
	\caption{Net displacement of the three-dimensional magneto-swimmer over one period as a function of the frequency.}

	\label{Fig:magneto-swimmer3d:res}
\end{figure}

The \cref{Fig:magneto-swimmer3d:res} shows the net displacement $\Delta x$ of the three-dimensional magneto-swimmer over one period as a function of the frequency $f$ of the external magnetic field. It leads to results similar to those in the two-dimensional case: the displacement increases with frequency, reaching a maximum at $f = 0.8$Hz, and then decreases. The behavior is consistent with the two-dimensional simulations and experimental observations reported in the literature \cite{oulmas_3d_2017}. Moreover, \cref{fig:intro:magneto} illustrates the deformation and motion of the micro-swimmer at different times during a stroke cycle. The transparent gray surfaces represent the boundaries between domains used by the processors in the parallel computation. However, the three-dimensional simulations are computationally 
expensive, requiring long execution times. One perspective of this work is to improve the solver and preconditioner strategies, as well as to 
verify the parallel scalability of the model.

\section{Perspectives and Conclusion}
\label{conclusion}

This work presents a comprehensive numerical framework for simulating fluid-structure interactions of elastic magneto-swimmers in confined domains, representing a significant advance toward developing digital twins for biomedical micro-robotics applications. The developed platform combines realistic modeling, numerical robustness, and high-performance computing capabilities to enable accurate simulation of magneto-swimmers.

The proposed full-order approach based on the Arbitrary Lagrangian-Eulerian formulation provides a major advancement in the physical understanding of active swimmer problems. The 2D and 3D validation benchmarks demonstrate the maturity and accuracy of the method, showing excellent agreement with experimental data from the literature. The framework successfully captures the complex nonlinear coupling between magnetic actuation, elastic deformations, and fluid interactions that characterize these systems.

Key contributions of this work include: (i) a robust finite element discretization strategy that handles large swimmer deformations through adaptive remeshing procedures, (ii) advanced coupling schemes that ensure stability and convergence of the fluid-structure interaction, and (iii) a parallel implementation within the open-source Feel++ library that leverages high-performance computing resources.

The developed platform opens numerous perspectives for biomedical applications. 
It incorporates advanced path-planning algorithms that allow autonomous navigation of magneto-swimmers through complex biological environments. These algorithms are essential for practical biomedical applications, where precise control and collision avoidance are critical for successful targeted delivery missions.
The framework also provides a solid foundation for in silico calibration using experimental data, optimal control strategies for targeted navigation, and the development of model reduction techniques. Its capabilities in handling realistic biological environments make it particularly well-suited for applications in drug delivery and micro-robot-assisted therapy.

Future developments will focus on extending the framework to multi-swimmer configurations to study collective behavior, implementing advanced control algorithms based on machine learning techniques, and optimizing the computational strategies for exascale computing. The integration with experimental validation campaigns will further enhance the platform's predictive capabilities, ultimately enabling the transition from laboratory prototypes to clinical applications in targeted biomedical interventions.

\section*{Acknowledgements}
This work of the Interdisciplinary Thematic Institute IRMIA++, as part of the ITI 2021-2028 program of the University of Strasbourg, CNRS and Inserm, was supported by IdEx Unistra (ANR-10-IDEX-0002), and by SFRI-STRAT’US project (ANR-20-SFRI-0012) under the framework of the French Investments for the Future Program.
The authors acknowledge the financial support of the French Agence Nationale de la Recherche (grant
ANR-21-CE45-0013 project NEMO), and Cemosis. Part of this work was also funded by the France 2030 NumPEx Exa-MA (ANR-22-EXNU-0002) project managed by the French National Research Agency (ANR).

\bibliographystyle{abbrv}
\bibliography{references}

\end{document}